\documentclass[final,1p,times]{elsarticle}
\usepackage{amssymb,amsmath}
\usepackage{float}
\newcommand\vp{\boldsymbol v'}
\newcommand\Sbb{\mathbb S}
\newcommand{\un}{u^{(n)}}

\newcommand{\fn}{f^{(n)}}

\newcommand{\hun}{\hu^{(n)}}
\newcommand{\hfn}{\hf^{(n)}}

\usepackage{amsmath}
\usepackage{amssymb}
\usepackage{amsthm}
\usepackage{graphicx}
\usepackage{psfrag}
\usepackage{fancyvrb}
\usepackage{color}
\usepackage{floatflt}
\usepackage{cancel}
\usepackage{nextpage}

\newcommand{\erf}{\text{erf}}

\definecolor{AndrasGreen}{rgb}{0, 0.7, 0}

\newcommand{\bv}{\boldsymbol{v}}

\newcommand{\bx}{\boldsymbol{x}}

\newcommand{\bz}{\boldsymbol{z}}

\newcommand{\bB}{\boldsymbol{B}}

\newcommand{\bE}{\boldsymbol{E}}

\newcommand{\hf}{\hat{f}}

\newcommand{\hu}{\hat{u}}

\def\barint_#1{\mathchoice
{\mathop{\vrule width 6pt height 3 pt depth -2.5pt
\kern -8.8pt \intop}\nolimits_{#1}}%
{\mathop{\vrule width 5pt height 3 pt depth -2.6pt
\kern -6.5pt \intop}\nolimits_{#1}}%
{\mathop{\vrule width 5pt height 3 pt depth -2.6pt
\kern -6pt \intop}\nolimits_{#1}}%
{\mathop{\vrule width 5pt height 3 pt depth -2.6pt
\kern -6pt \intop}\nolimits_{#1}}}

\newtheorem{prop}{Proposition}[section]
\newtheorem{cor}[prop]{Corollary}
\newtheorem{defn}[prop]{Definition}

\newcommand{\BeginAndrasProof}{\begin{proof}[{\bf Proof:}] \mbox{} \\[4pt]}
\newcommand{\EndAndrasProof}{\end{proof}}

\usepackage{graphicx}
\newtheorem{theorem}{Theorem}
\newtheorem{remark}{Remark}

\begin{document}
\begin{frontmatter}
\title{Fast elliptic solvers in cylindrical coordinates and the
Coulomb collision operator}
\author{Andras Pataki and Leslie Greengard}
\address{Courant Institute, New York University, 251 Mercer Street, 
NY, NY 10012}

\begin{abstract}
In this paper, we describe a new class of fast solvers for separable 
elliptic partial differential equations in cylindrical coordinates 
$(r,\theta,z)$ with free-space radiation conditions. 
By combining integral equation methods in the
radial variable $r$ with Fourier methods in $\theta$ and $z$, we show that
high-order accuracy can be achieved in both the governing potential
and its derivatives.  A weak singularity 
arises in the Fourier transform with respect to $z$ that is handled 
with special purpose quadratures.
We show how these solvers can be applied to the evaluation of the 
Coulomb collision operator in kinetic models of ionized gases.

\end{abstract}


\end{frontmatter}

\section{Introduction}

A variety of problems in computational physics require the 
solution of the Poisson and biharmonic equations in cylindrical
coordinate systems, particularly when the source distribution 
(the right-hand side) is axisymmetric or involves only a few
azimuthal modes. The present paper was motivated by the need to 
compute the Coulomb collision operator $C(f_a,f_b)$ in kinetic
simulations of the Boltzmann-Fokker-Planck equation
\cite{Bellan,BirdsallLangdon,ProcassiniBirdsall,Rosenbluth57,Shkarofsky}:
\begin{equation}
\partial_t f_a + \bv \cdot \nabla f_{a} +
\frac{e_a}{m_a} (\bE + \bv \times \bB) \cdot \partial_{\bv} f_a =
\sum_b C(f_a,f_b) .
\label{eq:boltzmann}
\end{equation}
Here, $f_a(\bx,\bv,t)$ denotes the state of an ionized gas for 
plasma species $a$ and the index $b$ runs over all species present. 
In the Fokker-Planck-Landau formalism \cite{Landau}, 
\begin{equation}  C(f_a,f_b) = \gamma_{ab} \partial_{\bv} \cdot \int 
\Sbb(\bv - \vp) \left( \frac{\partial_{\bv} f_a(\bv)}{m_a} f_b(\vp) -
f_a(\bv) \frac{\partial_{\vp} f_b(\vp)}{m_b} \right) \, d\vp
\end{equation}
where 
\begin{equation}
 \Sbb(\bv - \vp)_{ij} =  \delta_{ij} \frac{1}{| \bv - \vp |} -
\frac{ (v_i - v'_i) (v_j - v'_j)}{|\bv -\vp|^3} .
\end{equation}

An alternative representation makes use of the Rosenbluth potentials
\cite{Rosenbluth57}:

\begin{align}
	C(f_a,f_b) = \frac{\gamma_{ab}}{m_a} \partial_{\bv} \cdot 
        \left[ \partial_{\bv} \cdot (f_a \partial_{\bv}\partial_{\bv} G_b)
          - 2 \left(1 + \frac{m_a}{m_b}\right) f_a \partial_{\bv} H_b  \right]
        \label{eq:CRdef}
\end{align}
where
\begin{align}
	H_b(\bv) = \int \frac{1}{|\bv-\vp|} f_b(\vp) \, d\vp 
        \qquad {\rm or} \qquad \Delta H_b = - 4 \pi f_b   
        \label{eq:Hbdef}
\end{align}
and
\begin{align}
	G_b(\bv) = \int |\bv-\vp| f_b(\vp) \, d\vp  
	\qquad {\rm or} \qquad  \Delta^2 G_b = - 8 \pi f_b
	\label{eq:Gbdef}
\end{align}

Note that four derivatives of $G_b$ are required 
in (\ref{eq:CRdef}), while
$G_b$ itself satisfies the inhomogeneous biharmonic equation
(\ref{eq:Gbdef}).
Thus, direct discretization of the partial differential equation, followed by evaluation
of the collision operator via (\ref{eq:CRdef}) would require eight 
steps of numerical differentiation, with significant loss of accuracy.

It is natural, therefore, to consider alternative methods with the dual goals of 
achieving high order accuracy and minimizing the condition number of the
solution process. Because of the design of magnetic confinement devices
for plasmas, it is also important to be able to construct numerical 
methods in cylindrical coordinate systems, since the distribution
functions $f_b(\bv)$ are often axisymmetric or involve only a few 
azimuthal modes. 

There is, of course, a substantial literature on computing Coulomb
collisions and on solving elliptic partial differential equations in
cylindrical coordinates. We refer the reader to 
\cite{BirdsallLangdon,Caflisch,Coster,Hinton08,Karney,KilleenKerbelMcCoy,Lemou,Louche,
ProcassiniBirdsall,Rosenbluth57,Shkarofsky} 
for some methods in current use in plasma physics.
For a discussion of relativistic effects, see
\cite{Braams}.
Most closely related to our approach are the methods of 
\cite{FilbetPareschi,PareschiRusso} and \cite{JiHeld,KhabibrakhmanovKhazanov,KilleenKerbelMcCoy,Rosenbluth57}. The
first two are fast and achieve high order (``spectral") accuracy,
but use Fourier methods in Cartesian coordinates and do not address
the axisymmetric (or low azimuthal mode) case.
The latter rely on separation
of variables in spherical coordinates, for which the axisymmetric case leads naturally to a representation
involving Legendre polynomials and the general case to a representation involving
associated Legendre functions.

In the numerical analysis literature, most solvers based on
cylindrical coordinates tend to concern themselves with
periodic (in $z$) or finite domain boundary conditions rather than
free-space boundary conditions (see, for example \cite{Chen00,Lai07}).
Here, we develop a method for computing the Rosenbluth potentials 
using separation of variables and a mix of integral equation and 
Fourier analysis techniques. We show that free-space (radiation)
conditions can be imposed in a straightforward manner and that
high order accuracy can be achieved in all derivatives with minimal
loss of precision. The solver requires $O(N \log N)$ work, where
$N$ is the number of grid points used to sample the distribution
function.

Finally, we should make a remark about notation. The collision operator and the 
Rosenbluth potentials in (\ref{eq:Hbdef}),(\ref{eq:Gbdef}) 
are defined in velocity
variables, for which we will use the standard cylindrical coordinates $(r,\theta,z)$ for
$\bv$. In the context of plasma physics, 
 $r = |v_{\perp}|$, where $|v_{\perp}|$ is the magnitude
 of the component of the velocity perpendicular to the magnetic field, $\theta$ is the 
gyrophase angle, and $z = v_{||}$ is the component of the velocity field parallel to the magnetic field. The problem is purely axisymmetric when the velocity field is independent of the gyrophase angle. 

One disadvantage of our solver is that we can be adaptive in the $r$ direction,
but not in the $z$ or $\theta$ directions, since we use spectral discretizations in the latter
variables.
For fully adaptive three-dimensional calculations, one could employ fast multipole-accelerated integral equation solvers, as described in  \cite{Ethridge,Harper}. These methods directly
compute the convolution of the data $f_b(\bv)$ with the free-space Green's function.
In the axisymmetric case, one could use an axisymmetric version of the
fast multipole method \cite{Strickland}.
The constant, however, is larger for these schemes than for methods based on separation of variables, and we limit our attention to methods that rely on
a tensor product mesh in $r$, $\theta$ and $z$, which is adequate for most current
simulations of the Boltzmann-Fokker-Planck equation (\ref{eq:boltzmann}).

\section{The Poisson equation in cylindrical coordinates}

In order to compute the Rosenbluth potential $H_b$, we must solve the
Poisson equation in free space
\[
	\Delta u(\bv) = f(\bv).
\]
In cylindrical coordinates $\bv = (r, \theta, z)$,
we have
\begin{equation}
	u_{rr}(r,\theta,z) + \frac{1}{r} u_r(r,\theta,z)
	+ \frac{1}{r^2} u_{\theta \theta}(r,\theta,z) + u_{zz}(r,\theta,z) = f(r,\theta, z),
\label{eq:poisscyl}
\end{equation}
and we assume that $f$ is identically zero outside the region 
\[
	\Omega = \{(r,\theta,z): 0 \leq r \leq R, \;\; -A \leq z \leq A, \;\; 0 \leq \theta \leq 2\pi \}.
\]
Since $u$ and $f$ are periodic in $\theta$, we represent them as Fourier series:
\begin{align}
	u(r,\theta,z) &= \sum_{n=-\infty}^\infty \un(r,z) e^{in\theta}
        \\
	f(r,\theta,z) &= \sum_{n=-\infty}^\infty \fn(r,z) e^{in\theta}
	\label{eq:fmodes}
\end{align}

The derivatives in this representation will be written as 
\begin{align*}
	u_r(r,\theta,z) &= \sum_{n=-\infty}^\infty \un_r(r,z) e^{in\theta}
        &
	u_{rr}(r,\theta,z) &= \sum_{n=-\infty}^\infty \un_{rr}(r,z) e^{in\theta}
        \\
	u_{zz}(r,\theta,z) &= \sum_{n=-\infty}^\infty \un_{zz}(r,z) e^{in\theta}
        &
	u_{\theta \theta}(r,\theta,z) &= \sum_{n=-\infty}^\infty (-n^2) \un(r,z) e^{in\theta}
\end{align*}
Substituting into (\ref{eq:poisscyl}) and equating terms corresponding to the $n$th
azimuthal mode, we obtain:
\[
	\un_{rr}(r,z) + \frac{1}{r} \un_r(r,z) - \frac{n^2}{r^2} \un + \un_{zz}(r,z) = \fn(r, z).
\]
For each mode, we now have a partial differential equation (PDE)  
in the two variables $r$ and $z$ which we need
to solve on the rectangular domain 
\[ \Omega_{rz} = \{(r,z): 0 \leq r \leq R, \;\; -A \leq z \leq A \}. \]

Let us now take the Fourier transform of the equation in the $z$ direction,
That is we write
\begin{align}
	\un(r,z) &= \frac{1}{2 \pi} \int_{-\infty}^\infty e^{i\kappa z} \hun(r,\kappa) d\kappa 
	\label{eq:utrans} \\
	\hun(r,\kappa) &= \int_{-\infty}^\infty e^{-i\kappa z} \un(r,z) dz  \nonumber \\
	\fn(r,z) &= \frac{1}{2 \pi} \int_{-\infty}^\infty e^{i\kappa z} \hfn(r,\kappa) d\kappa 
	\label{eq:ftrans} \\
	\hfn(r,\kappa) &= \int_{-\infty}^\infty e^{-i\kappa z} \fn(r,z) dz \nonumber
\end{align}

In the Fourier transform domain, the PDE becomes an ordinary differential 
equation (ODE), where $\kappa$ (as well as $n$) is now fixed:
\begin{align}
	\hun_{rr}(r,\kappa) + \frac{1}{r} \hun_r(r,\kappa) - \left( \frac{n^2}{r^2}
        + \kappa^2 \right) \hun(r,\kappa) &= \hfn(r,\kappa).
\label{eq:modbessel}
\end{align}
To simplify notation (when the context is clear),
we will write $\hu(r)$ instead of $\hun(r,\kappa)$ and $\hu'(r)$ instead of 
$\hun_r(r,\kappa)$ to denote the derivative when discussing the solution
of the ODE. 

The equation (\ref{eq:modbessel}) is an inhomogeneous modified Bessel equation 
\cite{AS}. In the homogeneous case, the equation has
two linearly independent solutions, namely $I_n(|\kappa|r)$ and $K_n(|\kappa|r)$,
the modified Bessel functions of order $n$.
The function $I_n(|\kappa|r)$ is regular at the origin, and grows exponentially as
$r \rightarrow \infty$, while 
$K_n(|\kappa|r)$ is logarithmically singular at the origin, but decays exponentially fast
as $r \rightarrow \infty$.

\subsection{Boundary conditions for the modified Bessel equation} \label{sec:bcbessel}

In order to have a properly posed ODE, we seek two boundary conditions,
one at $r=0$ and one at $r=R$, beyond which the equation is homogeneous.
For the $n=0$ mode, the condition
 \[ \hu'(0) =  \hu_r^{(0)}(0,\kappa) = 0 \]
ensures regularity at the origin, while for modes $n \neq 0$
 \[ \hu(0) = \hun(0,\kappa) = 0 \]
is necessary. This is easily seen from taking the limit of the equation
(\ref{eq:modbessel}) as $r \rightarrow 0$.

Since we are seeking to solve the Poisson equation in free space, our ODE is actually 
posed on the half line $[0,\infty]$, with the radiation condition
that the solution decay at infinity. This can be accounted for exactly in terms of a
suitable boundary condition at $r=R$. To see this, note that for $r>R$ the solution must be proportional to
$K_n(|\kappa|r)$, since $I_n(\kappa|r)$ grows without bound. That is, 
\[
	\hu(r) =  C_{n,\kappa} \cdot K_n(|\kappa|r)  \hspace{1in}
	{\rm for}\ r  \geq R,
\]
where $C_{n,\kappa}$ is an unknown constant. The solution on $[0,R]$ and its
derivative must match this solution at $r=R$, so that
\begin{align*}
 \hu(R) &= C_{n,\kappa} \cdot K_n(|\kappa|R),	\\
 \hu'(R) &= C_{n,\kappa} |\kappa| \cdot K_n'(|\kappa|R).
\end{align*}
Eliminating the constant $C_{n,\kappa}$
we obtain the exact ``radiation" boundary condition:
\begin{equation}
	\hu(R) - \frac{K_n(|\kappa|R)}{|\kappa| \cdot K_n'(|\kappa|R)} \hu'(R) = 0.
\end{equation}
In summary, the ODE boundary value problem we must solve is
(\ref{eq:modbessel}), subject to the boundary conditions:
\begin{align}
        \hu'(0) & = 0  && n = 0 \nonumber
        \\
        \hu(0) & = 0  && n \neq 0
	\\
	\hu(R) - \frac{K_n(|\kappa|R)}{|\kappa| \cdot K_n'(|\kappa|R)} \hu'(R) & = 0.
	\label{eq:nrbc}
\end{align}

In broad terms, this completes the description of the Poisson solver, which proceeds in 
four steps.

\vspace{.2in}

\hrule

\vspace{.05in}

\centerline{\sf Informal description of fast Poisson solver}

\vspace{.05in}

\hrule

\vspace{.05in}

\begin{enumerate}
\vspace{-4pt}
\setlength{\itemsep}{2pt}
\setlength{\parskip}{0pt}
\item Expand the right hand size $f(r,\theta,z)$ as a Fourier series in $\theta$, in
order to get $\fn(r,z)$,
\item Compute the Fourier transform of $\fn(r,z)$ in the $z$ direction to get $\hfn(r,\kappa)$,
\item Solve the ODE (\ref{eq:modbessel}) for each $\kappa$ and $n$ to obtain
$\hun(r,\kappa)$,
\item Compute the inverse Fourier transform of $\hun(r,\kappa)$ to get $\un(r,z)$,
\item Sum the Fourier series in $\theta$ to get the final solution $u(r, \theta, z)$.

\vspace{.1in}

\hrule
\end{enumerate}

\vspace{.1in}

We will rely on fairly standard methods for all of the above, except Steps 3 and 4.
For Step 3, we use an analytic solution based on knowledge of the underlying Green's
function for the ODE and accelerated by a simple ``sweeping" algorithm. Step 4
will require some care, since it is straightforward to show that
$\hun(r,\kappa)$ is logarithmically
singular as $\kappa \rightarrow 0$ for $n=0$ and has a singularity of the 
order $\kappa^{2|n|} \log  \kappa$ for $n \neq 0$. 

\section{Discretization and solution} \label{sec:discretization}

\vspace{.1in}

We assume $f(r,\theta,z)$ is given on a tensor product grid with 
$N_\theta$ equispaced points in the $\theta$ direction on $[0,2\pi]$, $N_z$ 
equispaced points in the $z$ direction on $[-A,A]$, 
and $N_r$ points in the $r$ direction on $[0,R]$. We divide 
$[0,R]$ into $N_I$ intervals with interval endpoints $R_0 = 0, R_1, R_2, \dots, R_{N_I} = R$.
We use a $P$th order (scaled) Chebyshev grid on each, so that
$N_r = N_I \, P$. We will denote by $\{ r_j \, | j = 1,\dot,N_r\}$ the grid points in increasing order.
When the particular interval $m$ ($1 \leq m \leq N_I$) is of interest, the $p$th grid
point on that interval $(1 \leq p \leq P$) is 
$r_j = r_{ (m-1) \, P + p}$.

The discretized data will be denoted by 
\[
	f_h(r_j, \theta_n, z_k) = f(r_j, \theta_n, z_k) \ {\rm for}\ 
          0 \leq j < N_r,\  0 \leq n < N_\theta,\   0 \leq k < N_z.
\]

\subsection{Step 1: Transformation in $\theta$} \label{sec:step1}

\vspace{.1in}

We use the fast Fourier transform (FFT) to compute $\fn_h(r_j, z_k)$,
the discretized version of $\fn(r,z)$:
\begin{equation}
	\fn_h(r_j, z_k) = \frac{2\pi}{N_\theta} \, 
	\sum_{l=0}^{N_\theta - 1} e^{-\frac{2\pi i}{N_\theta} n l} f_h(r_j, \theta_l, z_k)
	\approx \fn(r_j,z_k).
\label{eq:gntrap}
\end{equation}

It should be noted that, if $f(r,\theta,z)$ is $n$-times differentiable, then the 
series (\ref{eq:fmodes}) truncated after $N$ terms has an error of the order
\begin{equation}
 O \left(\frac{1}{N^{n-1}}\right) \qquad\hbox{as}\quad N \rightarrow \infty\ .
\end{equation}
If $f$ is infinitely differentiable, then the error goes to zero faster than
any finite power of $1/N$. Schemes with this property are often referred to as having 
{\em spectral accuracy}. Moreover, the trapezoldal rule approximations of the series coefficients
in (\ref{eq:gntrap}) converge at the same rate \cite{GO,Trefethen}.

\subsection{Step 2: Transformation in $z$} \label{sec:step2}

\vspace{.1in}

Since $f(r, \theta, z)$ and  $\fn(r, z)$ are compactly supported, we need to compute the
finite integral
\[
	\hfn(r,\kappa) = \int_{-A}^A e^{-i\kappa z} \fn(r,z) dz .
\]
Letting  $h_z = \frac{2A}{N_z}$ and $z_l = lh_z$, the trapezoidal rule yields:
\begin{equation}
\hfn_h \left(r_j,\frac{\pi}{A} k\right) =
	\frac{2A}{N_z} \sum_{l=-N_z/2}^{N_z/2-1} e^{-\frac{2\pi i}{N_z} l k} \fn_h(r_j,z_l)
	\approx \hfn \left(r_j,  \frac{\pi}{A} k\ \right).
\label{eq:ztrap}
\end{equation}
This is computable using the FFT, and yields the values of the Fourier transform at 
equally spaced points of step size $\pi/A$ in the $\kappa$ domain.
A few remarks are in order:

\begin{itemize}
\item The ratio $\frac{N_z}{A}$ determines the range of frequencies that are resolved.
If $\frac{N_z}{A}$ increases ($h_z$ decreases),
higher frequency modes of the data are computed. 

\item
We will assume that, to precision $\epsilon$,
$\hf(r,\kappa)$ is supported on the interval $[-\kappa_{\max}, \kappa_{\max}]$.
For a given $A$, $N_z$ must be chosen sufficiently large that
$\pi N_z/(2A) > \kappa_{\max}$. (This is simply asking that the grid in $z$ be fine enough
to resolve the data.)
\item
Increasing $N_z$ and $A$ simultaneously so that $N_z/A$ remains fixed
leaves the range of $\kappa$ unchanged, but increases the number of sample points
where $\hfn_h$ is computed in the range $[-\pi N_z/(2A), \pi N_z/(2A)]$. 
\item
The trapezoidal approximation (\ref{eq:ztrap}) is spectrally accurate, since the
integrand and all its derivatives are assumed to have vanished by the time $z = \pm A$.
\end{itemize}

\subsection{Step 3: Solving the modified Bessel equation} \label{sec:ODEsolve}

\vspace{.1in}

We turn now to the solution of the modified Bessel equation
(\ref{eq:modbessel}), subject to the boundary conditions (\ref{eq:nrbc})
for $\kappa \neq 0$. (As noted above, the equation has a weakly singular solution
at $\kappa=0$. Our quadrature rule for computing the inverse Fourier transform  
in section \ref{sec:invFT} will avoid the origin when integrating along the $\kappa$ axis.)
\\
One possible approach to solving the equation is to use a spectral integration-based  ODE 
solver \cite{spectral1}
that represents the second derivative as a Chebyshev series:
\[
	\hu''(r) = \sum_{k=0}^N \alpha_k T_k(r) \, .
\]
Multiplying the equation (\ref{eq:modbessel}) by $r^2$ and 
systematic use of the following two 
identities for Chebyshev polynomials
\[
	\int T_n(r) dr + C = \frac{1}{2(n+1)} T_{n+1}(r) - \frac{1}{2(n-1)} T_{n-1}(r)
\]
\[
        r T_n(r) = \frac{T_{n+1}(r) + T_{n-1}(r)}{2}
\]
yields a banded linear system (of bandwidth 7) to which are appended two dense rows that correspond to 
the imposition of the desired boundary conditions.
Such a system can be solved in linear time by careful Gaussian elimination, achieving spectral
accuracy.
For non-singular ODEs, this linear system can be viewed as the discretization of a second-kind 
integral equation for the unknown second derivative, and thus as a well-conditioned formulation of the 
problem. Unfortunately, in our case, the differential operator is singular at the origin. As a result, 
the integral equation is not of the second kind and the approach becomes ill-conditioned for fine grids,
with the attendant loss of precision.
\\[6pt]
An alternative strategy is to use the fact that our ODE is classical and well studied, with a 
known Green's function $G^n_\kappa(r,s)$. We can, therefore,
write down the exact solution as a convolution:
\begin{equation}
\hun(r,\kappa) = \int_0^R  G^n_\kappa(r,s) f(s) ds 
\label{ode2exact}
\end{equation}
where
\begin{align*}
 G^n_\kappa(r,s) &=   \left\{ \begin{array}{lll} 
              I_n(\kappa r)K_n(\kappa s)/W(s) & \mbox{if} & r \leq s \\[4pt]
              K_n(\kappa r)I_n(\kappa s)/W(s) & \mbox{if} & s < r 
          \end{array}\right.  \ {\rm where} \\
           W(s) &= \kappa ( I_n'(\kappa s)K_n(\kappa s) - K_n'(\kappa s)I_n(\kappa s))
           = - \frac{1}{s}
\end{align*}
This choice of Green's function correctly imposes the regularity condition at the origin
and the radiation condition at infinity. In this formulation, there is no need to solve a linear
system - one needs only to evaluate the integral in (\ref{ode2exact}).
Naive implementation of this formula would require $O(N_r^2)$ work.
Because of the structure of the Green's function, however, there is a simple $O(N_r)$ solver
based on the observation that

\begin{align}
\hun(r,\kappa) =   & K_n(\kappa r) \int_0^r  I_n(\kappa s) f(s)/W(s) ds \;+  \nonumber \\ 
& I_n(\kappa r) \int_r^R  K_n(\kappa s) f(s)/W(s) ds \, .
\label{eq:sweep}
\end{align}

The only source of error comes from the quadrature approximation of the preceding integrals.
Derivatives of the solution are also obtained analytically. For example,

\begin{align}
\hun_r(r,\kappa) =  & \kappa K_n'(\kappa r) \int_0^r  I_n(\kappa s) f(s)/W(s) ds \;+ \nonumber  \\
& \kappa I_n'(\kappa r) \int_r^R  K_n(\kappa s) f(s)/W(s) ds \, .
\label{eq:sweepd}
\end{align}

There are some implementation issues in using (\ref{eq:sweep}), having to do with 
scaling and quadrature
due to the fast growth/decay of Bessel functions for increasing $n$ and $r$.
In particular, when $r$ lies in the $m$th interval denoted by 
$[R_{m-1},R_{m}]$, we write
\begin{align}
\hun(r,\kappa) =  & \frac{K_n(\kappa r)}{K_n(\kappa R_m)} \,
 \int_0^r  I_n(\kappa s) K_n(\kappa R_m) \, f(s)/W(s) ds \; +  \nonumber \\
& I_n(\kappa r) K_n(\kappa R_m) \int_r^R  \frac{K_n(\kappa s)}{K_n(\kappa R_m)} \, 
f(s)/W(s) ds \, .
\label{eq:sweepscale}
\end{align}

\subsection{Step 4: Computing the inverse Fourier transform} \label{sec:invFT}

We now need to compute the inverse Fourier transform of $\hun(r,\kappa)$
to recover $\un(r,z)$, according to (\ref{eq:utrans}).
Since $\hun(r,\kappa)$ is compactly supported to the desired precision on $[-\kappa_{\max}, \kappa_{\max}]$, we actually need to compute
\begin{align}
	\un(r,z) \approx \frac{1}{2\pi} \int_{-\frac{\pi}{2A}N_z}^{\frac{\pi}{2A} N_z} 
	e^{i\kappa z} \, \hun(r,\kappa) d\kappa \,
	\label{eqn-Axi-IFT-1}
\end{align}
where (as discussed in section \ref{sec:step2}) $\pi N_z/(2A) > \kappa_{max}$.
A complication is that $\hun(r,\kappa)$ has a logarithmic singularity at $\kappa=0$.

Fortunately, in the last decade or so, a variety of quadrature rules have been developed
that rely on slight modifications of the trapezoidal rule, yield high-order accuracy, and
still permit the use of the FFT. Two such schemes are the end-point corrected
trapezoidal rule due to Kapur and Rokhlin \cite{KR} and the hybrid Gauss-trapezoidal rule
due to Alpert \cite{Alpert}. We will make use of the latter.

\begin{theorem} (modified from \cite{Alpert}). Let $f(\kappa)$ be a compactly supported
function on $[-\kappa_{max},\kappa_{max}]$ which is smooth away from the origin and 
takes the form 
\[  f(\kappa) = s_1(\kappa) log(|\kappa|) + s_2(\kappa)  \]
in a neighborhood of the origin, where $s_1$ and $s_2$ are smooth functions. 
Let 
\[ I(f) = \int_{-\kappa_{max}}^{\kappa_{max}}  f(\kappa) \, d\kappa \]
and let $h = \frac{2 \kappa_{\max}}{N_z}$.
Then, for every integer $m>0$ and every $N_z > 2m$, there exist weights $w_{l,m}$ and nodes $\kappa_{l,m}$
such that 
\begin{equation}
 I_h(f) = {h} \, \sum_{\substack{k=-N_z/2 \\ |k| \ge m}}^{N_z/2} 
 f(kh) + 
\sum_{\substack{l=-m}}^{m} 
 w_{l,m} f(\kappa_{l,m})
 \label{alpertrule}
\end{equation}
satisfies 
\[ I_h(f) = I(f) + O(h^{m}) \, . \]
\end{theorem}

In other words, the hybrid Gauss-trapezoidal rule achieves $m^{th}$ order accuracy
by replacing the $2m$ trapezoidal nodes nearest the origin with specially located
nodes (and weights). The paper \cite{Alpert} provides tables of these nodes for orders
2-16 (and the corresponding ones for a variety of other singularities as well).

In the present context, therefore, we will compute 
the integral (\ref{eqn-Axi-IFT-1}) using the formula (\ref{alpertrule}):
\begin{equation}
\un(r,z)
	= \frac{h}{2\pi} 
	\sum_{\substack{k=-N_z/2 \\ |k| \ge m}}^{N_z/2} 
e^{i\kappa_k z} \hun(r,\kappa_k)  + 
\sum_{\substack{l=-m}}^{m} 
 w_{l,m}  e^{i\kappa_{l,m} z} \hun(r,\kappa_{l,m}) 
\label{kappa_discr}
\end{equation}
with mesh spacing $h = 2\pi N_z/(2A) \div N_z$, $\kappa_k =  \frac{\pi}{A} k$, and 
$w_{l,m}, \kappa_{l,m}$ taken from \cite{Alpert}.
Evaluating  $\un(r,z)$ on our grid $z_j = \frac{2A}{N_z} j$, we have:
\[
\un(r,z_j)
	= \frac{1}{2A} 
		\sum_{\substack{k=-N_z/2 \\ |k| \ge m}}^{N_z/2} 
	 e^{\frac{2\pi i \, jk}{N_z}} \hun(r,\kappa_k)
	 +
	 \sum_{\substack{l=-m}}^{m} 
 w_{l,m}  e^{i\kappa_{l,m} \frac{2Aj}{N_z} } \hun(r,\kappa_{l,m}) 
\]
The first term is straightforward to compute with the (inverse) FFT, requiring $N_z \log N_z$
operations. The second term can be computed directly using 
$O(N_z m)$ operations, where $m$ is the order of the quadrature rule.
For $m$ sufficiently large, the sums can be computed simultaneously using the 
non-uniform FFT (see \cite{dutt} and the more recent review \cite{nufft1}).

\begin{remark}
The quadrature rule (\ref{kappa_discr}) determines the discrete values of 
the continuous  Fourier transform variable $\kappa$ where $\hun(r,\kappa)$
needs to be sampled. The number of such points is $O(N_z+m)$. This, in turn,
tells us where $\hfn(r,\kappa)$ is needed. The values at the 
regular nodes are obtained with the FFT, as discussed in section \ref{sec:step2}. 
The values at the irregular nodes $\kappa_{l,m}$ can be computed directly or 
using the non-uniform FFT.
\end{remark}

\begin{remark} {\bf (Oversampling)}\, . \ 
In practice, there is one more issue which needs to be addressed. 
In the integral (\ref{eqn-Axi-IFT-1}), $z$ is bounded by $A$, so that the 
most oscillatory integrand is $e^{i\kappa A} \hun(r,\kappa)$. It is easy to see that
there are a maximum of $N_z/2$ periods of the function $e^{i\kappa A}$ 
over the interval of integration $[-\pi N_z/(2A),\pi N_z/(2A)]$. The maximum
for the function $\hun(r,\kappa)$ is similar.
Thus,  the trapezoidal rule with $N_z$ points yields only one point per wavelength
for the most oscillatory argument, in violation of the Shannon sampling theorem.
We, therefore, {\em oversample} the integrand by a factor
of $\eta$, by setting $N_z' = \eta N_z$. (As discussed in 
section \ref{sec:step2}, we must simultaneously set
 $A' = \eta A$ in computing the forward transform.)

Setting $\eta=1$ yields exponentially small errors near z=0, but $O(1)$ errors at $z = A$.
Setting $\eta \geq 2$  ensures convergence for  $z$  in the 
entire range $[-A,A]$, with exponential improvement as $\eta$ increases.
Setting $\eta=4$ is sufficient for double precision accuracy for  $N_z > 16$,
assuming the function is bandlimited to machine precision at $\kappa_{max}
= \pi N_z/(2A)$. 
\end{remark}

\subsection{Step 5: Sum the Fourier in $\theta$ to obtain the full solution}
 \label{sec:sumFS}

This is completely straightforward. 
As in Step 1, we may use the fast Fourier transform (FFT) to compute $u(r_j, \theta,z_k)$
at equispaced points $\theta_l = \frac{2 \pi \, l}{N_\theta}$:
\begin{equation}
	u(r_j,\theta_l,z_k) \approx \frac{1}{N_\theta} \, 
	\sum_{n=0}^{N_\theta - 1} e^{\frac{2\pi i}{N_\theta} n l} \un(r_j, z_k) \, .
\label{eq:sumFS}
\end{equation}

\subsection{Computing derivatives of the solution}

One useful feature of spectral solvers is that derivative are straightforward
to compute with high order accuracy.

\begin{enumerate}
\item {\em First and second $r$-derivatives}:
Our ODE solver returns both the solution $\hun(r,\kappa)$
and its derivatives $\hun_{r}(r,\kappa)$ and $\hun_{rr}(r,\kappa)$ on our grid.
Thus, we can compute $u_r(r,\theta,z)$
and $u_{rr}(r,\theta,z)$ by the same technique as for $u(r,\theta,z)$:
evaluating the inverse $z$-Fourier transform
and the $\theta$ Fourier series for $\hun_{r}$ and $\hun_{rr}$, respectively.

\item {\em $\bz$-derivatives}:
In the present paper, 
$z$ derivatives  are obtained through multiplication by $i\kappa$ in the inverse Fourier transform step:
\[
	\un(r,z) = \frac{1}{2\pi} \int_{-\infty}^\infty e^{i\kappa z} \hun(r,\kappa) d\kappa
        \qquad\Rightarrow \qquad
	\frac{\partial^m}{\partial z^m} \un(r,z)
        = \frac{1}{2\pi} \int_{-\infty}^\infty (i\kappa)^m e^{i\kappa z} \hun(r,\kappa) d\kappa
\]
and the quadrature rule described above for logarithmic singularities.

\item {\em $\theta$-derivatives}:
In the present paper, we also compute $\theta$-derivatives spectrally, by differentiating
the Fourier series:
\[
	u(r,\theta,z) = \sum_{n=-\infty}^\infty \un(r,z) e^{in\theta}
        \qquad \Rightarrow \qquad
	\frac{\partial^m}{\partial z^m} u(r,\theta,z) = \sum_{n=-\infty}^\infty (in)^m \un(r,z) e^{in\theta},
\]
using the FFT.

\item {\em Mixed derivatives}:
Since derivatives in $r, \theta, z$ are computed at independent steps of the algorithm,
they are easily combined.
For example, if we want to compute $u_{rzz}$, we start with 
$\hun_r(r, \kappa)$ and compute the inverse $z$ Fourier transform on the function
$-\kappa^2 \hun_r(r,\kappa)$ to get $\un_{rzz}(r,z)$, followed by evaluating the
Fourier series in the $\theta$ direction via the FFT.
\end{enumerate}

\begin{remark}
One can easily obtain
$z$ derivatives without numerical differentiation,
once $u$, $u_r$, $u_{rr}$ and $u_{\theta \theta}$ are known. The original
PDE (\ref{eq:poisscyl}) becomes a second order ODE in $z$, and the method
of spectral integration \cite{spectral1} can be applied directly.
We have not implemented this option, since the condition number of 
Fourier differentiation is only $O(N)$, so that with $1000$ points in $z$ (or
$1000$ azimuthal modes), one can still obtain at least 10 digits of accuracy in double 
precision.
\end{remark}

\section{The biharmonic equation in cylindrical coordinates}

For the Rosenbluth potential $G_b$, we must solve the
biharmonic equation in free space
\[
	\Delta^2 u(\bv) = f(\bv).
\]
In cylindrical coordinates $\bv = (r, \theta, z)$, after Fourier transformation in 
$z$ and $\theta$, we obtain the fourth order Bessel-type ODE:
\begin{align}
	u_{rrrr} + \frac{2}{r} \hun_{rrr}
        - \left( \frac{1 + 2n^2}{r^2} + 2\kappa^2 \right) \hun_{rr}
        + \left( \frac{1 + 2n^2}{r^3} - \frac{2\kappa^2}{r} \right) \hun_r
        \hspace{0.5in} \nonumber
        \\
        + \left( \frac{n^4 - 4n^2}{r^4} + \frac{2k^2 n^2}{r^2} + \kappa^4 \right) \hun
        &= \hfn
        \label{eq:biharmcyl}
\end{align}

\subsection{Boundary conditions for the fourth order Bessel-type equation} \label{sec:bcfourth}

Since (\ref{eq:biharmcyl}) is fourth order,
we need four boundary conditions to have a properly posed ODE. We impose
two at $r=0$ and two at $r=R$, beyond which the equation is homogeneous.
To ensure regularity at the origin, 
for the $n=0$ mode, it is sufficient to impose
 \[ \hu_r^{(0)}(0,\kappa) = \hu_{rrr}^{(0)}(0,\kappa) = 0. \]
For the $n=1$ mode, we set
 \[ \hu^{(1)}(0,\kappa) = \hu_{rr}^{(1)}(0,\kappa) = 0 \, , \]
and for modes $n \geq 2$, we set
 \[ \hu^{(n)}(0,\kappa) = \hu_{r}^{(n)}(0,\kappa) = 0  \, .\]

These conditions are easily derived by taking the limit of the equation
(\ref{eq:modbessel}) as $r \rightarrow 0$ and the fact that the null space of the 
differential operator is spanned by 
\[ \{ I_n(|\kappa|r), r \, I_n'(|\kappa|r),  K_n(|\kappa|r),  r \, K_n'(|\kappa|r)  \}.
\]
It remains to determine a radiation condition at $r=R$, so that the derivative of the 
solution is bounded at infinity. As for the Poisson equation, we proceed 
by observing that the solution $\hu^{(n)}(r)$ for $r>R$ must take the form
\[
	\hu(r) = C^1_{n,\kappa}  K_n(|\kappa|r) + C^2_{n,\kappa}  r K_n'(|\kappa|r) \, ,
\]
where $C^1_{n,\kappa}, C^2_{n,\kappa}$ are unknown constants.
This follows since the derivatives of $I_n(\kappa|r)$ and $r  I_n'(\kappa|r)$ 
grow without bound. 
The solution on $[0,R]$ and its
derivative must match this solution at $r=R$, so that
\begin{align*}
 \hu^{(n)}(R) &= C^1_{n,\kappa} \cdot K_n(|\kappa|R) +
 C^2_{n,\kappa} \cdot r K_n'(|\kappa|R),	\\
 \hu^{(n)}_r(R) &= C^1_{n,\kappa} \cdot \kappa K_n'(|\kappa|R) +
 C^2_{n,\kappa} \cdot (K_n'(|\kappa|R) +  \kappa r K_n''(|\kappa|R) ), \\
  \hu^{(n)}_{rr}(R) &= C^1_{n,\kappa} \cdot \kappa^2 K_n''(|\kappa|R) +
 C^2_{n,\kappa} \cdot (2 \kappa K_n''(|\kappa|R) +  \kappa^2 r K_n'''(|\kappa|R) ), \\
  \hu^{(n)}_{rrr}(R) &= C^1_{n,\kappa} \cdot \kappa^3 K_n'''(|\kappa|R) +
 C^2_{n,\kappa} \cdot (3 \kappa^2 K_n'''(|\kappa|R) +  \kappa^3 r K_n''''(|\kappa|R) ).
\end{align*}
Eliminating the constants, we obtain two exact ``radiation" boundary conditions to be imposed
on the combination of $\hun(R)$ and its first three derivatives. The formula is complex and 
omitted, since we won't use it. We will instead use an exact solution based on the Green's function. 

\subsection{Discretization and solution}
The solution of the biharmonic equation is analogous to that of the Poisson equation,
so we just highlight the differences.
\\[6pt]
After separation of variables, we need to solve a fourth order Bessel type equation.
As before, we could proceed by expanding the highest derivative in a Chebyshev series and integrating,
but the resulting linear system again loses precision because of the singular nature of the 
differential operator at the origin. (The loss is, in fact, much more severe than for the second order
(Poisson) equation.)
\\[6pt]
Alternatively, we can construct the Green's function for the ODE using the linearly independent
fundamental solutions
$I_n(\kappa r),\, rI_n'(\kappa r),\, K_n(\kappa r),\, r K_n'(\kappa r)$, imposing the 
regularity condition at $r=0$ and the decay condition as $r\rightarrow \infty$:
\begin{align*}
	G^{(n)}_\kappa(r,s)
	&= \left\{ \begin{array}{ll}
            \left[ I_n(\kappa r) s K_n'(\kappa s) + r I_n'(\kappa r) K_n(\kappa s) \right] / W(s),
            &r \leq s
            \\[4pt]
            \left[ K_n(\kappa r) s I_n'(\kappa s) + r K_n'(\kappa r) I_n(\kappa s) \right] / W(s),
            &r > s
            \end{array} \right.
          \hspace{24pt} \text{with} \hspace{24pt}
          W(s) = - \frac{2 \kappa}{s}
\end{align*}
The solution involves computing four integrals (instead of two). The sweeping 
method is virtually the same as that used for the Poisson equation.
\begin{align*}
	\hun(r)
	&= K_n(\kappa r) \int_0^r s I_n'(\kappa s) f(s) / W(s) ds
          \;+\; r K_n'(\kappa r) \int_0^r I_n(\kappa s) f(s) / W(s) ds
          \\
          &\hspace{0.5in}
          \;+\; I_n(\kappa r) \int_r^\infty s K_n'(\kappa s) f(s) / W(s) ds
          \;+\; r I_n'(\kappa r) \int_r^\infty K_n(\kappa s) f(s) / W(s) ds \, .
\end{align*}
The Fourier transform of the solution in $z$ has a more severe singularity in the biharmonic case,
due to the fact that the free-space Green's function does not decay. Fortunately, however,
we are only interested in second derivatives of the biharmonic potential, and they have only logarithmic
singularities, so our special-purpose quadratures from section \ref{sec:invFT}
yield the desired accuracy. (More elaborate methods involving singularity subtraction could be developed
if one wanted the biharmonic potential or its first derivatives.)

\section{The collision operator}
\label{sec:collisionop}
Now that we've described how to solve the Poisson and biharmonic equations (\ref{eq:Hbdef})
and (\ref{eq:Gbdef}), we turn our attention to the collision operator (\ref{eq:CRdef}).
If we express all the derivative terms in cylindrical coordinates,
the axisymmetric collision operator becomes:
\begin{align}
	C(f^a, f^b) &= \frac{\gamma_{ab}}{m_a} \left[ C_b(f^a, f^b) - 2 \left( 1 + \frac{m_a}{m_b} \right) C_p(f^a, f^b) \right]
        \label{eqn:collisioncyl}
        \\
        C_p(f^a, f^b)
        &   
	= - 4 \pi f^a f^b + f^a_r H^b_r + f^a_z H^b_z
        \nonumber
        \\
        C_b(f^a, f^b)
        &= - 8 \pi f^a f^b + f^a_r \left[ 2 G^b_{rzz} + 2 G^b_{rrr} + \frac{2}{r} G^b_{rr} - \frac{1}{r^2} G^b_r \right]
        + f^a_z \left[ \frac{2}{r} G^b_{rz} + 2 G^b_{rrz} + 2G^b_{zzz} \right]
        \\
        & \hspace{3in}
        + f^a_{rr} G^b_{rr} + 2 f^a_{rz} G^b_{rz} + f^a_{zz} G^b_{zz}
        \nonumber
\end{align}

\section{Numerical Examples}
\label{sec:numex}

In order to test the convergence of the algorithm, it is desirable to compare the results to a nontrivial
exact solution. For a right-hand side consisting of a radially symmetric Gaussian:
\begin{align*}
	f(\rho) &= \frac{E}{(4 \pi v)^{3/2}}
	&
        E &= e^{-\frac{\rho^2}{4v}}, \;\;
        \rho^2 = x^2 + y^2 + z^2  \, .
\end{align*}
we can compute the exact solution
to both the Poisson and the biharmonic equations
\[
	\Delta u = f  \qquad  \Delta^2 v = f     
\]
as well as to the components of  the collision operator $C_p,C_b$:
\begin{align}
	  u(\rho)	&= - \frac{R}{4 \pi \rho} 
	&
        R = \erf\left(\frac{\rho}{2\sqrt{v}} \right), 
                \nonumber \\
        v(\rho)
        &= - \left[\frac{\rho}{8\pi} + \frac{v}{4 \pi \rho} \right] R
        - \frac{\sqrt{v}}{4 \pi^{3/2}} E
        \nonumber
        \\
	C_p &= - \frac{E^2}{8 \pi^2 v^3} + \frac{ER}{16 \pi^{3/2} v^{5/2} \rho},
        &
	C_b = - \frac{E^2}{2 \pi^2 v^3} + \frac{ER}{4 \pi^{3/2} v^{5/2} \rho}.
        \label{eqn:gausscoll}
\end{align}
After a change of variables to cylindrical coordinates,
we can find explicit analytic formulas
for all quantities produced by our solvers
(although some of the formulas need to be treated  carefully to avoid catastrophic
cancellations in their numerical evaluation).

\begin{figure}[H]
\begin{center}
\includegraphics[width=5.0in,angle=0]{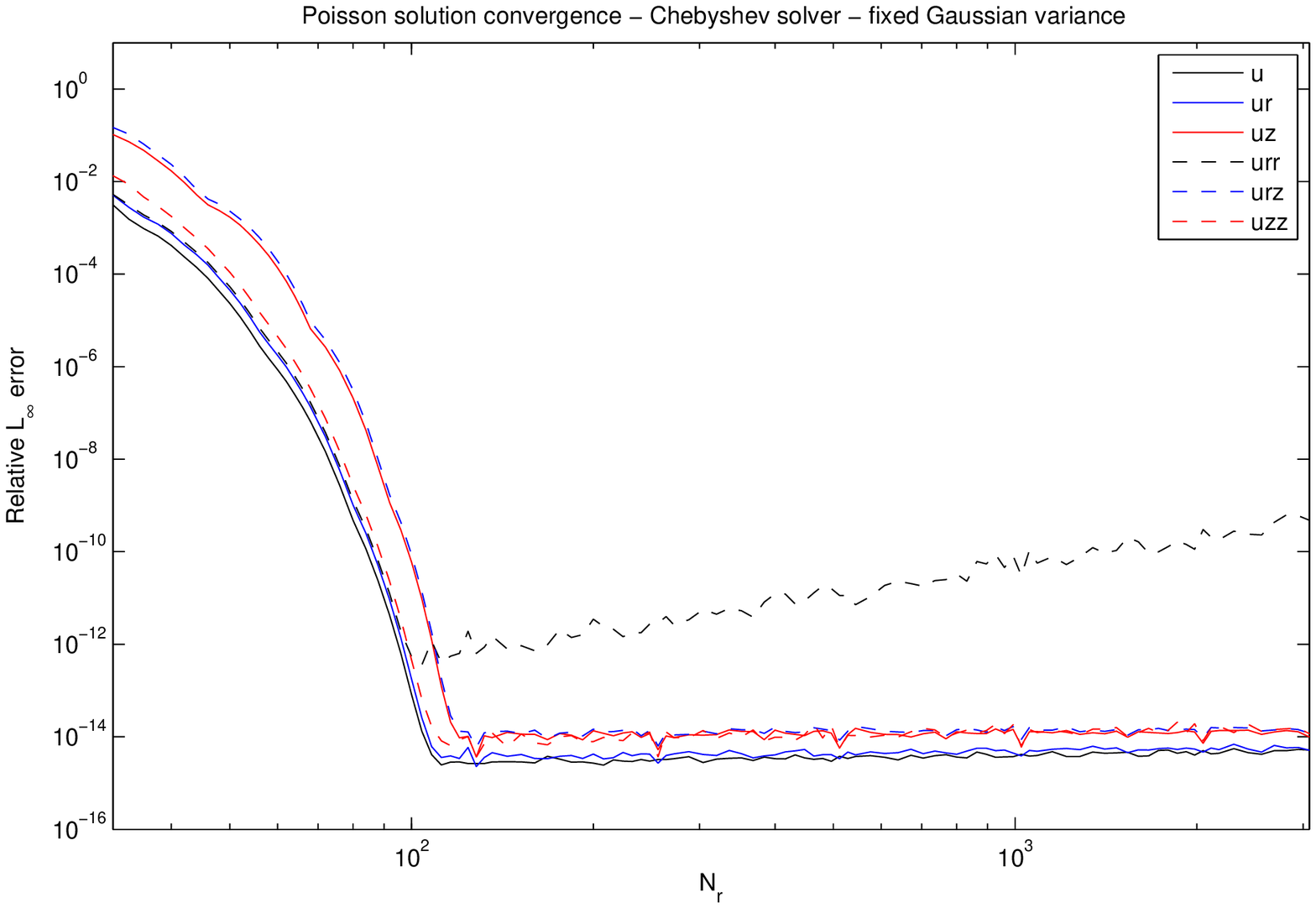}
\includegraphics[width=5.0in,angle=0]{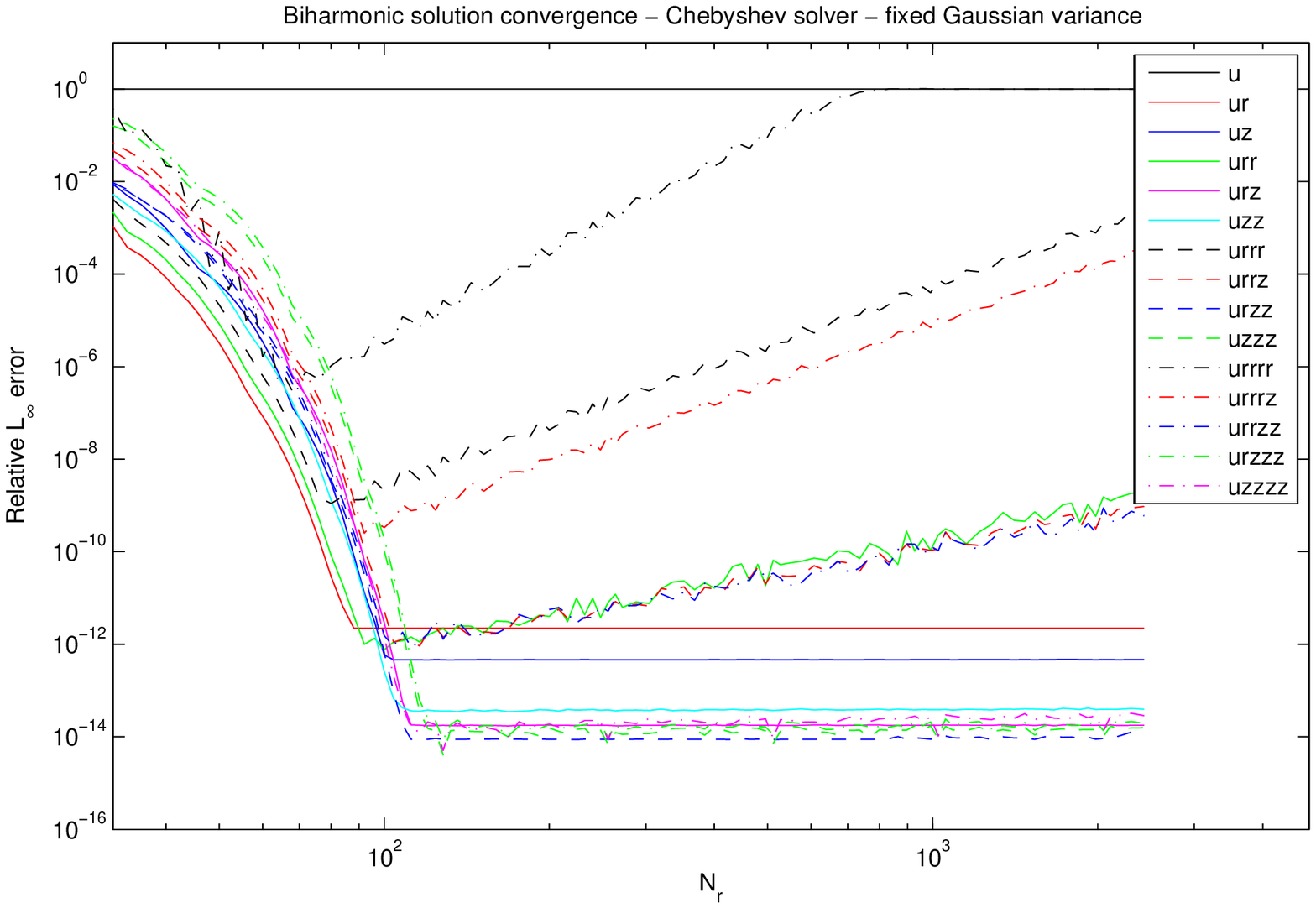}
\end{center}
\caption{The relative $L_\infty$ errors for the spectral integration based solvers, when applied
to the Poisson and biharmonic equations (top and bottom, respectively).} \label{fig:1}
\end{figure}

\clearpage
{\bf Example 1:} Let us first consider the convergence of the solver for a single Gaussian
of variance $v = 0.223$.
Using the Chebyshev (spectral integration) solvers, note that after an initially rapid convergence,
the higher derivatives start to diverge due to the ill-conditioning of the linear system (Fig \ref{fig:1}).

\begin{figure}[H]
\begin{center}
\includegraphics[width=5.0in,angle=0]{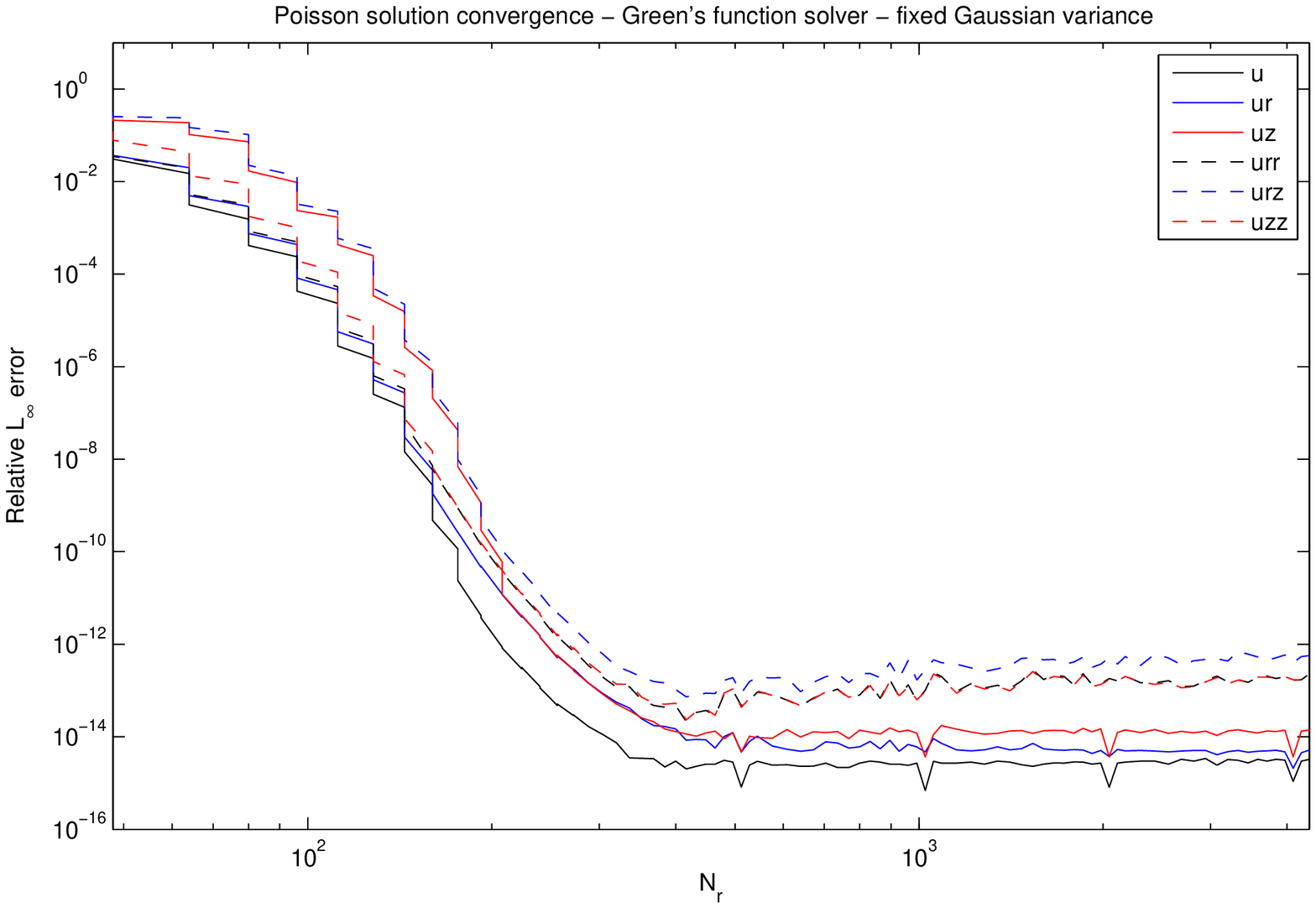}
\includegraphics[width=5.0in,angle=0]{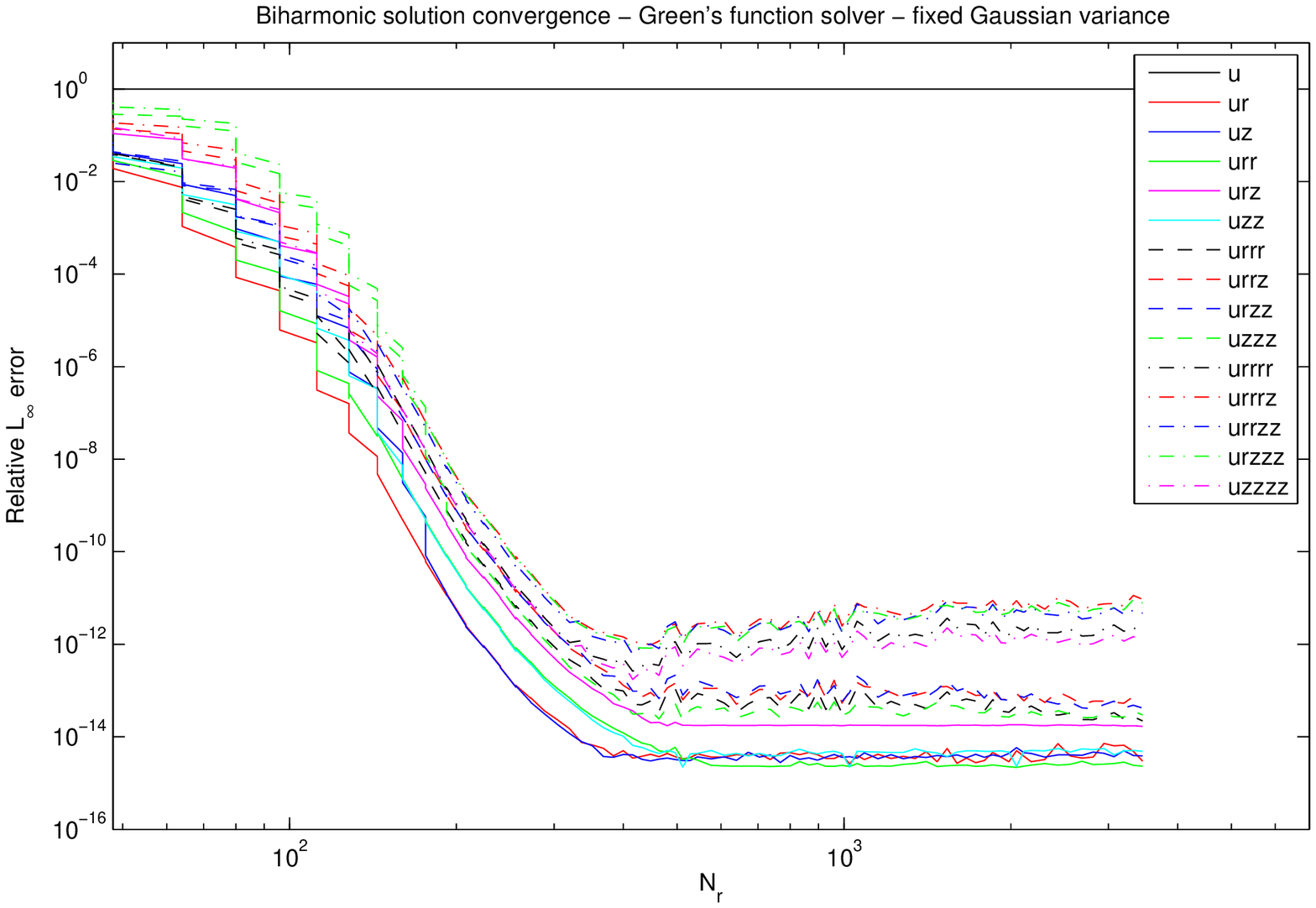}
\end{center}
\caption{The relative $L_\infty$ errors for the Green's function based solvers, when applied
to the Poisson and biharmonic equations (top and bottom, respectively).} \label{fig:2}
\end{figure}

\clearpage
Using the Green's function based ODE solver, we see
that the the solution and its derivatives converge spectrally without the loss of precision for fine grids 
that affects the spectral integration based schemes (Fig \ref{fig:2}).

\begin{figure}[H]
\begin{center}
\includegraphics[width=5.0in,angle=0]{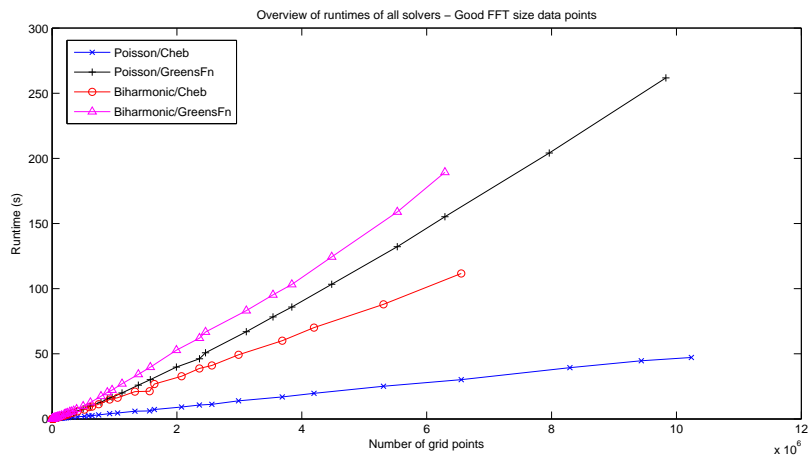}
\end{center}
\caption{Runtime for the Green's function and spectral integration based solvers, when applied
to the Poisson and biharmonic equations.} \label{fig:3}
\end{figure}

\clearpage
Fig. \ref{fig:3} depicts the runtime performance of the solver. In the present implementation,
the spectral integration based code for the Poisson equation requires about 25 seconds for
5 million grid points, while the biharmonic solver is about three times slower. This is within
a small factor of the performance of FFT-based codes for doubly or triply periodic constant
coefficient PDEs on regular grids. The Green's function based solvers are a bit slower
at present. We estimate that straightforward optimization/precomputation could result in a
factor of 2-3 speed-up.

\clearpage
{\bf Example 2:} 
The Coulomb collision operator appears as a source term in the Boltzmann-Fokker-Planck equation (\ref{eq:boltzmann}).
Not considering convective and electromagnetic terms, it affects the evolution of ionized gas
consisting of a single species as
\[
	\partial_t f + (\text{convective and electromagnetic terms})
        = C(f, f)
\]
For a single species of ion, the Maxwellian distribution (a Gaussian centered at the origin)
is an equilibrium state.
From the formulas (\ref{eqn:collisioncyl}) and  (\ref{eqn:gausscoll}),
it can be verified analytically that $C(f^G, f^G) = 0$.
Computing this result numerically on a $192\times 128$ grid we obtain zero to about 14 digits.
Note that $C_p$ and $C_b$ do not vanish independently. There is a real cancellation between
the two contributions when inserted into the formula (\ref{eqn:collisioncyl}).
\begin{center}
\begin{tabular}{|c|c|c|}
	\hline
        &&
        \\[-6pt]
	$\displaystyle ||C_p(f^G,f^G)||_\infty$
	& $\displaystyle ||C_b(f^G,f^G)||_\infty$
	& $\displaystyle ||C(f^G,f^G)||_\infty$
        \\[4pt] \hline
        && \\[-6pt]
        $3.403 \cdot 10^{-3}$ &  $1.361 \cdot 10^{-2}$ &  $5.851 \cdot 10^{-14}$
        \\[2pt] \hline
\end{tabular}
\end{center}

To illustrate the diffusive nature of the collision operator,
let us consider a perturbation to the equilibrium solution, by constructing
an anisotropic Gaussian density, with slightly different variances
in the $r$ and $z$ directions:
\[
	f(\rho) = \frac{1}{\sqrt{4 \pi v_z}\, 4 \pi v_r}  \, e^{ -(x^2+y^2)/(4 v_r)} \, e^{ -z^2/(4 v_z)} 
\]
where $v_r =  1.107$ and $v_z =  1.353$.

\begin{center}
\begin{figure}[H]
\includegraphics[width=5.5in,angle=0]{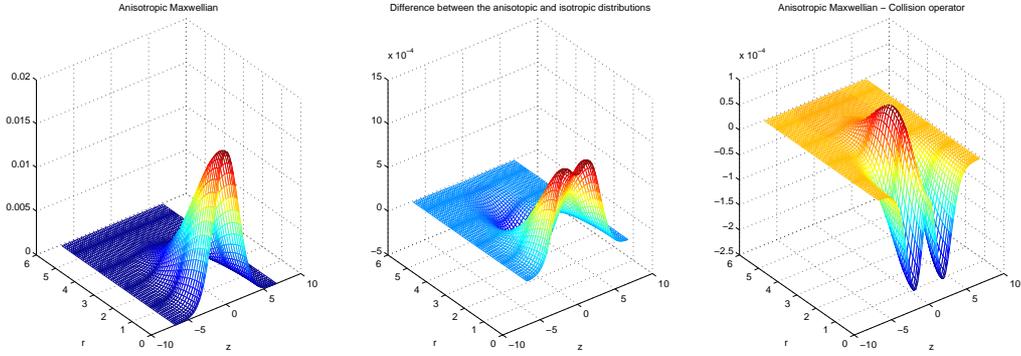}
\caption{Anisotropic Maxwellian, its difference from a Maxwellian and its computed collision operator}
\end{figure}
\end{center}
Notice that where the anisotropic distribution is too small (the valley in the central plot),
the collision operator is positive, thus it tends to increase $f$.
Where the anisotropic distribution function is too large (the peaks in the central plot),
the collision operator is negative, tending to decrease $f$.
Thus, the collision operator indeed has the effect of moving an anisotropic Maxwellian
towards an isotropic equilibrium distribution.

\clearpage
To demonstrate the full three dimensional solver, we placed three Gaussian source densities
at different locations in $R^3$, two with positive weight and one with negative weight.
The following plot shows the solution overlaid on the $r-\theta-z$ grid
as a three dimensional contour surface plot with an octant ``cut out'':

\begin{figure}[H]
\begin{center}
\includegraphics[width=5.0in,angle=0]{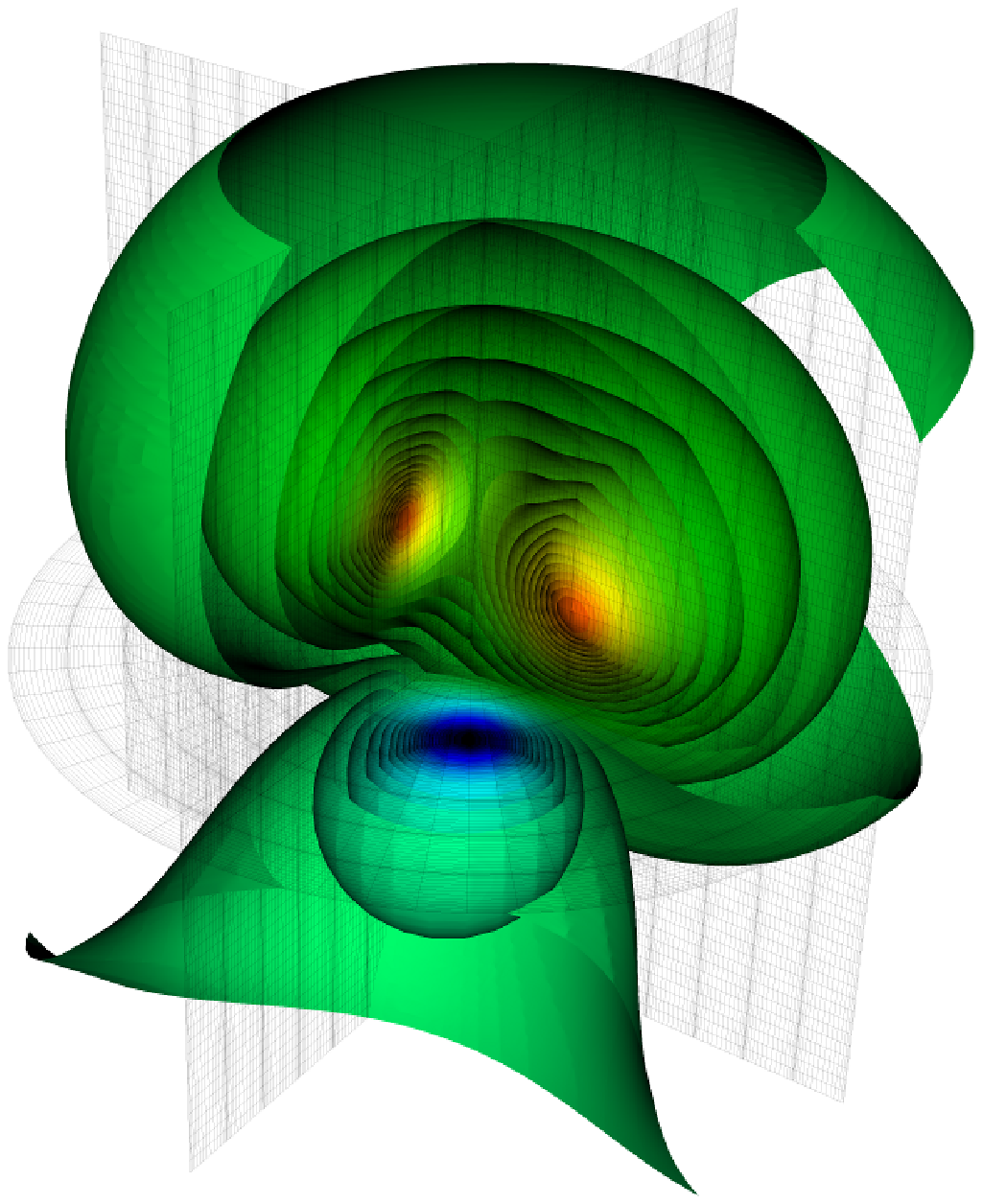}
\end{center}
\end{figure}

We used Gaussian variances of 0.2, 0.6 and 0.3 centered at the points with Cartesian coordinates
$(4.3,1.2,3.6)$, $(-1.1,4.1,-0.8)$, and $(5.3,3.5,-3.2)$ and weights $-1.0$, $-1.3$, and $1.2$, respectively. With a grid in $(r,\theta,z)$ of $192 \times 64 \times 96$ on $[0,16] \times [0, 2\pi] \times [-16,16]$, 
we obtained 12 digits of accuracy in the solution and its gradient. The total execution
time was 18 secs. on a single core of a 2.5GHZ CPU, using an oversampling factor in $z$ of 4.

\section{Conclusion}
\label{sec:conclusion}

This paper describes a new fast solver for separable 
elliptic partial differential equations in cylindrical coordinates that is both fast and high-order accurate,
with solution times comparable to a few FFTs using the same number of degrees of freedom.
Combined with the Rosenbluth formalism, it permits the rapid evaluation of the Coulomb collision operator
in kinetic models of ionized gases.

Our solver is particularly useful when the number of azimuthal modes is small. For full three-dimensional
problems, it is quite likely that Cartesian-based methods will be more effective, particularly since one
can use fast multipole-based, fully adaptive solvers. Here, we require regular grids in the $\theta$ and $z$
directions, which is sufficient for many applications. 

Several open problems remain. One involves the construction of fast, fully implicit collision operators, so that large time steps can be taken in the Boltzmann-Fokker-Planck  equation (for which there is already a significant literature). Another involves the development of  methods for solving elliptic partial differential equations in complicated axisymmetric geometries rather than in free space - that is, interior or exterior to a surface of revolution. These problems are currently being investigated, with progress to be reported at a later date.

\section*{Acknowledgements}
We thank C.S. Chang and Eisung Yoon for several useful discussions. 
This work was supported by 
the National Science Foundation under grant DMS06-02235 and the Department of Energy 
under contract DEFG0288ER25053.

\bibliographystyle{hplain}

\end{document}